\documentclass[11pt]{amsart}
\usepackage{rlepsf,graphicx,latexsym,amssymb}
\usepackage[all]{xy}
\CompileMatrices
\theoremstyle{plain}
\newtheorem{theorem}{Theorem}[section]
\newtheorem{lemma}[theorem]{Lemma}
\newtheorem{proposition}[theorem]{Proposition}
\newtheorem{corollary}[theorem]{Corollary}
\newtheorem{definition}[theorem]{Definition}
\theoremstyle{remark}
\newtheorem{example}[theorem]{Example}

\numberwithin{equation}{section}
\numberwithin{figure}{section}

\def \N{\mathbb{N}}
\def \Z{\mathbb{Z}}
\def \Q{\mathbb{Q}}
\def \F{\mathbb{F}}

\def \A{\mathcal{A}}
\def \Cyl{\mathcal{C}yl}
\def \L{\mathcal{L}}

\def \deg{{\rm deg}}

\def \incl{{\rm incl}}
\def \int{{\rm int}}
\def \Id{{\rm Id}}

\def \Ker{{\rm Ker}}
\def \mod{{\rm mod}}

\def \Tors{{\rm Tors}}
\def \P{{\rm P}}
\def \U{{\rm U}}
\def \V{{\rm V}}

\begin{document}

\title[]{Finite-type invariants of three-manifolds\\ and the dimension subgroup problem}

\date{November 13, 2006}

\author[]{Gw\'ena\"el Massuyeau}

\address{Institut de Recherche Math\'ematique Avanc\'ee,
CNRS/Universit\'e Louis Pasteur,  7 rue Ren\'e Descartes,
67084 Strasbourg, France}
\email{massuyeau@math.u-strasbg.fr}
\subjclass[2000]{57M27, 16S34}
\keywords{$3$-manifold, finite-type invariant, group ring, N-series, dimension subgroup}

\begin{abstract}
For a certain class of compact oriented $3$-manifolds, M. Goussarov and K. Habiro have conjectured
that the information carried by finite-type invariants should be characterized
in terms of ``cut-and-paste'' operations defined by the lower central series of the Torelli group of a surface.
In this paper, we observe that this is a variation of a classical
problem in group theory, namely the ``dimension subgroup problem.''
This viewpoint allows us to prove, by purely algebraic methods,  
an analogue of the Goussarov--Habiro conjecture for finite-type invariants with values in a fixed field.
We deduce that their original conjecture is true at least in a weaker form.
\end{abstract}

\maketitle
\tableofcontents

\section{Introduction}

The Goussarov--Habiro theory is aimed at understanding how $3$-manifolds can be obtained one from the other
by cut-and-paste operations of a certain kind \cite{Goussarov,Habiro,Goussarov_bis,GGP}.
In particular, it applies to the study of finite-type invariants introduced by Ohtsuki \cite{Ohtsuki}:
The latter are invariants of $3$-manifolds which, in a sense, 
behave polynomially with respect to certain surgery operations.

We will work with compact oriented $3$-manifolds whose boundary, if any,
is identified with a fixed abstract surface. Those manifolds are considered up
to homeomorphisms that preserve the orientation and the boundary identification.
The kind of surgery modifications that are used in the Goussarov--Habiro theory 
can be described as follows: Given a compact oriented $3$-manifold $M$,
a handlebody $H \subset M$ and a Torelli automorphism $h$ of $\partial H$
(that is, $h:\partial H \to \partial H$ is a homeomorphism
which acts trivially in homology), one can form a new compact oriented $3$-manifold:
$$
M_h := (M \setminus \int\ H) \cup_h H.
$$
The move $M \leadsto M_h$ is called a \emph{Torelli surgery}.

A possible characterization of a \emph{finite-type invariant} is as follows.
Let $f$ be an invariant of compact oriented $3$-manifolds with values in $A$, an Abelian group.
Then, $f$ is a finite-type invariant of degree at most $d$ if,
for any manifold $M$ and for any set $\Gamma$ of $d+1$ pairwise disjoint handlebodies in $M$ --
each coming with a Torelli automorphism of its boundary -- the following identity holds:
\begin{equation}
\label{eq:fti}
\sum_{\Gamma' \subset \Gamma} (-1)^{|\Gamma'|} \cdot f\left(M_{\Gamma'}\right) =0 \quad \in A.
\end{equation}
Here, $M_{\Gamma'}$ denotes the manifold obtained from $M$
by the simultaneous Torelli surgery along the handlebodies belonging to $\Gamma'$.

There are plenties of finite-type invariants \cite{LMO}
and the natural problem is to ``quantify'' how fine they are, degree by degree.
For this, Goussarov and Habiro have introduced, for every integer $k\geq 1$, 
the \emph{$Y_k$-equivalence}. This equivalence relation among manifolds can be characterized as follows:
Two compact oriented $3$-manifolds $M$ and $M'$ are $Y_k$-equivalent if $M'$ can be obtained from $M$ by a Torelli surgery
$$
M \leadsto M_h = M'
$$
such that $h$ belongs to the $k$-th term of the lower central series of the Torelli group of $\partial H$.\\[-0.3cm]

\noindent
\textbf{Fact.} \emph{If two manifolds are $Y_{d+1}$-equivalent, then
they are not distinguished by finite-type invariants of degree at most $d$.}\\[-0.3cm]

\noindent
The converse has been proved for integral homology $3$-spheres by Habiro and Goussarov \cite{Habiro, Goussarov_bis} 
which provides, in this special case, a geometric characterization of
the power of finite-type invariants. Unfortunately,
the converse does not hold in general.\footnote{We have inserted a counter-example below, in Section \ref{sec:GH}.}

Nevertheless, there are  manifolds of a special type
which play a key role in the Goussarov--Habiro theory.
Let $\Sigma$ be a connected compact oriented surface. Recall from \cite{Habiro, Goussarov} that
a \emph{homology cylinder} over $\Sigma$ is a cobordism $M$ from $-\Sigma$ to $\Sigma$
(with corners, if $\partial \Sigma \neq \varnothing$) 
which can be obtained from the cylinder $\Sigma\times[-1,1]$ by a Torelli surgery.
Homology cylinders can be composed, and their monoid is denoted by $\Cyl(\Sigma)$.\\[-0.3cm]

\noindent
\textbf{Conjecture} (Goussarov--Habiro)\textbf{.} 
\emph{Two homology cylinders over $\Sigma$ are $Y_{d+1}$-equivalent if, and only if,
finite-type invariants of degree at most $d$ do not separate them.}\\

In this paper, we identify the algebra underlying the ``Goussarov--Habiro conjecture'' (or GHC, for short)
by observing that it is a special instance of a general problem in group theory,
which is known as the ``dimension subgroup problem'' (or DSP, for short).
Classically, for a given group $G$, the question is to decide whether 
the $i$-th subgroup of the lower central series coincides with the $i$-th \emph{dimension subgroup},
i.e$.$ this subgroup determined by the $i$-th power of the augmentation ideal of $\Z[G]$.
Here, we consider the DSP in the more general setting where
the lower central series is replaced by an \emph{N-series} in Lazard's sense \cite{Lazard}.
See Section \ref{sec:DSP} for a precise statement of that problem.

In Section \ref{sec:GH}, the GHC is reduced to the DSP.
This is only an algebraic re-formulation of topological results, 
which have been proved independently by Habiro \cite{Habiro} and Goussarov \cite{Goussarov_bis} 
by means of their ``calculus of claspers'' and which notably led to their result on homology spheres.
In particular, they have shown that the monoid $\Cyl(\Sigma)$ 
quotiented out by the $Y_{d+1}$-equivalence relation is a group, 
and this is the group $G$ in which we are interested. 
Unfortunately, the DSP having no general answer, 
this algebraic reduction does not solve by itself the GHC (other than up to degree $d=2$).

However, the DSP has been solved for coefficients in a field.
In the case of the lower central series, this is due to Mal'cev, Jennings and Hall in the zero characteristic case
\cite{Mal'cev,Jennings_bis,Hall} 
and to Jennings and Lazard in the positive characteristic case \cite{Jennings,Lazard}.
Following Passi's book \cite{Passi}, these results are extended to any N-series in Section \ref{sec:weak}.
Thus, we obtain a version of the GHC for finite-type invariants with values in a given field $\F$;
the statement depends on the characteristic of $\F$:

\begin{theorem}
\label{th:rational_GH}
Let $\F$ be a field of characteristic $0$,
and let $M,M'$ be two homology cylinders over  $\Sigma$.
Finite-type invariants of degree at most $d$ with values in $\F$ do not distinguish $M$ from $M'$
if, and only if, there exists an $n\in \N^*$ such that $M^n$ is $Y_{d+1}$-equivalent to ${M'}^n$.
\end{theorem}

\begin{theorem}
\label{th:modular_GH}
Let $\F$ be a field of characteristic $p>0$,
and let $M,M'$ be two homology cylinders over  $\Sigma$.
Finite-type invariants of degree at most $d$ with values in $\F$ do not distinguish $M$ from $M'$
if, and only if, there exist some $C_1,\dots,C_r\in \Cyl(\Sigma)$ 
and some $e_1,\dots,e_r \in \N_0$ such that
\begin{itemize}
\item $C_i$ is $Y_{k_i}$-equivalent to $\Sigma \times [-1,1]$ for some  $k_i\in \N^*$
such that $k_i \cdot p^{e_i} \geq d+1$,
\item $M' = M \cdot \prod_{i=1}^r {C_i}^{p^{e_i}}$.
\end{itemize}
\end{theorem}

Next, a weak converse to the above Fact is deduced for homology cylinders:

\begin{corollary} 
\label{cor:big_degree}
There exists an integer $\mathbf{d}(\Sigma,d)\geq d$
such that, if  two homology cylinders over $\Sigma$ are not distinguished by finite-type invariants
of degree at most $\mathbf{d}(\Sigma,d)$, then they are $Y_{d+1}$-equivalent.
\end{corollary}

\noindent
This weak version of the GHC has the following consequence:
Two homology cylinders are $Y_k$-equivalent for all $k\geq 1$ if, and only if, 
they are not distinguished by finite-type invariants.

In Section \ref{sec:Quillen}, Theorem \ref{th:rational_GH} and Theorem \ref{th:modular_GH}
are extended as follows: The graded algebra dual to $\F$-valued finite-type invariants 
is the enveloping algebra of the Lie algebra of homology cylinders 
introduced by Habiro in \cite{Habiro}. To show this, we must extend 
the result proved by Quillen for the lower central series in \cite{Quillen} to an arbitrary N-series.

Last section is an appendix where the above characterizations of a finite-type invariant
and of the $Y_k$-equivalence relation, are shown to agree with the original definitions 
introduced by Goussarov and Habiro. But, this is well-known.\\

Before going into the proofs, we add a few more comments about the GHC. Recall that
one of the advantages of the Goussarov--Habiro theory is to apply
to links as well as manifolds and to support
-- at the level of results and techniques  -- 
a strong analogy between the two cases of study.
Thus, there exist for string links exact analogues of the above results, 
and these are proved following exactly the same algebraic method. 

Actually, the conjecture had been stated by Habiro for string links  \cite[Conjecture 6.13]{Habiro},
and had been announced by Goussarov as a result for homology cylinders \cite[Theorem 4]{Goussarov}.
But, as far as we know, no proof has never been communicated and
the question has been recently asked by Polyak in \cite[Conjecture 10.17]{Ohtsuki_problems}.

In another paper, Goussarov proved a weak version 
of the GHC for string links \cite[Theorem 10.4]{Goussarov_bis}, 
which is quite different from our Corollary \ref{cor:big_degree}
in the sense that it uses a notion of \emph{partially defined} finite-type invariant.
Besides, at the  end of the same paper \cite[Remark 10.8]{Goussarov_bis}, 
he announced a result that seems to be very close to Corollary \ref{cor:big_degree}.\\

\noindent
\emph{Acknowledgments.} We are sincerely grateful to Kazuo Habiro for stimulating and helpful conversations.
We also thank Jean--Baptiste Meilhan for his comments on the paper.
This research has been carried out while we were visiting the RIMS (Kyoto University) with a JSPS fellowship.

\section{The dimension subgroup problem}

\label{sec:DSP}

In this section, we briefly review the dimension subgroup problem.
The reader is refered to the survey \cite[\S 2]{Hartley} for an overview of that problem, 
and to the books \cite{Passi,Passman} for details and proofs.

Let $G$ be a group.  

\begin{definition}
An \emph{N-series} $N= (N_i)_{i\geq 1}$ of $G$
is a descending chain of subgroups
$$
G=N_1 \geq N_2 \geq N_3 \geq \cdots
$$
such that
$$
\left[N_i,N_j\right] \subset N_{i+j}  \quad \quad \forall i,j\geq 1.
$$
Moreover, if $G/N_i$ is torsion-free for all $i$, then the N-series $N$ is said to be \emph{$0$-restricted}.
If, for some fixed $p\in \N^*$, $x\in N_i$ implies that $x^p \in N_{ip}$  for all $i$,
then $N$ is \emph{$p$-restricted}.
\end{definition}

\noindent
Observe that an $N$-series is always central. 

\begin{example}
The lower central series $\gamma(G)$ of $G$, defined by
$\gamma_1(G):=G$ and $\gamma_{i+1}(G) := [\gamma_i(G),G]$, is an N-series.
\end{example}

Let $R$ be a commutative ring.
The augmentation ideal of the group ring $R[G]$ is denoted by $I(R[G])$.

\begin{definition}
A \emph{filtration} $\Delta = (\Delta_i)_{i\geq 1}$ of $I\left(R[G]\right)$
is a descending chain of $R$-submodules
$$
I\left(R[G]\right) = \Delta_1 \geq \Delta_2 \geq \Delta_3 \geq \cdots
$$
such that
$$
\Delta_i \cdot \Delta_j \subset \Delta_{i+j} \quad \quad \forall i,j\geq 1.
$$
\end{definition}

\noindent
Note that each term $\Delta_i$ of a filtration $\Delta$ is a two-sided ideal of $R[G]$.  

\begin{lemma}
\label{lem:filtrations}
If $\Delta$ is a filtration of $I\left(R[G]\right)$, then the sequence 
$\left(G\cap (1 + \Delta_i)\right)_{i\geq 1}$ is an $N$-series of $G$. 
Moreover, if $R$ is a field of characteristic $p\geq 0$, 
then it is $p$-restricted.
\end{lemma}

\begin{proof} Those are basic facts \cite{Passi,Passman}. 
Yet, we recall the identities on which the proof is based, since they will be used later.
Set $D_i:= G\cap (1 + \Delta_i)$. 

First of all, for every $g,h \in G$, one has that
$$
g h^{-1} -1 = \left( (g-1) - ( h - 1)\right) \cdot h^{-1} 
$$
which proves that $D_i$ is a subgroup of $G$.

Second, for all $g,h \in G$, one has that
\begin{equation}
\label{eq:comm}
[g,h] - 1 = g^{-1} h^{-1} \left((g-1)(h-1) - (h-1)(g-1)\right) \quad \in R[G]
\end{equation}
which proves that $[D_i,D_j] \subset D_{i+j}$.

Third, for all $g\in G$ and assuming that $R$ is a field of characteristic $p>0$, one has
\begin{equation}
\label{eq:binom_p}
g^p - 1 = (g-1)^p \quad \in R[G]
\end{equation}
which proves the implication $\left(g\in D_i \Rightarrow g^p \in D_{ip}\right)$.

Finally, for all $g\in G$, $m\in \N^*$ and provided that $R$ is a field of characteristic $0$, one has that
\begin{equation}
\label{eq:binom_0}
g-1 = \frac{1}{m}(g^m-1) - \frac{1}{m}\sum_{k=2}^m \binom{m}{k} (g-1)^k \quad \in R[G]
\end{equation}
which shows, after an induction on $i\geq 1$, that $G/D_i$ is torsion-free.
\end{proof}

Conversely, an N-series $N$ of $G$ 
induces a filtration $\Delta(N;R)$ of $I\left(R[G]\right)$ defined by
$$
\Delta_i(N;R) =  \left\langle\ \left(g_1-1\right) \cdots \left(g_r-1\right) \
\left| \ g_1 \in N_{k_1},\dots, g_r \in N_{k_r}\ \hbox{ and }\
 k_1+ \cdots + k_r \geq i \right.\ \right\rangle_R.
$$
Next, by Lemma \ref{lem:filtrations},
the filtration $\Delta(N;R)$ induces itself an N-series $D(N;R)$ of $G$ defined by
$$
D_i(N;R) := G\cap \left(1 + \Delta_i(N;R)\right).
$$

\begin{definition}
Let $N$ be an N-series of $G$. The \emph{$i$-th dimension subgroup}
of $G$ defined by $N$  with coefficients in $R$ is $D_i(N;R)$.
\end{definition}

\noindent
\textbf{Dimension Subgroup Problem.}
\emph{Compute the N-series $D(N;R)$ in terms of the original N-series
$N$ and the coefficient ring $R$}.\\[-0.3cm]

\noindent
Note the obvious inclusion
$N_i \subset D_i(N;R)$,
but the converse does not hold in general:

\begin{example}
\label{ex:classical_DSP}
Classically, the DSP refers to the special case
where $N = \gamma(G)$ is the lower central series of $G$, in which case one has  
$$
\Delta_i(\gamma(G);R) = I(R[G])^i
$$ 
as follows from (\ref{eq:comm}).
The most difficult problem is when $R=\Z$ is the ring of integers,  
and it had been conjectured for thirty years that
\begin{equation}
\label{eq:DSP}
\gamma_i(G) \stackrel{?}{=} G \cap \left( 1+ I(\Z[G])^i \right).
\end{equation}
The first counter-example is due to Rips (for a finite $2$-group $G$ and $i=4$, see \cite{Rips}).
\end{example}

Nevertheless, in some cases, the two N-series $N$ and $D(N;R)$
may happen to coincide, as illustrated by the following results.

\begin{lemma}
\label{lem:Abelian}
Assume that $G$ is Abelian and let $N$ be an N-series of $G$.
Then, one has that $D_i(N;\Z) = N_i$ for all $i\geq 1$.
\end{lemma}

\begin{proof}
Let us denote additively the law of the Abelian group $G/N_i$. The canonical projection 
$G \to G/N_i$ extends additively to a map $\pi_i: \Z[G] \to G/N_i$, for which one observes that
$$
\pi_i((x-1)(y-1)) = \pi_i(xy) - \pi_i(x) - \pi_i(y) + \pi_i(1) =
(\pi_i(x) + \pi_i(y)) - \pi_i(x) - \pi_i(y) + 0=0.
$$
So, $\pi_i$ vanishes on $\Delta_i(N;\Z)$. We conclude that, for all $g\in D_i(N;\Z)$,
$\pi_i(g-1)=\pi_i(g) -\pi_i(1)=\pi_i(g)$ vanishes or, equivalently, that $g\in N_i$.
\end{proof}

In general, it is true that $D_i(N;\Z)=N_i$ for $i=1,2$. This is trivial for $i=1$ 
and follows from Lemma \ref{lem:Abelian} for $i=2$ (by reducing to $G/N_2$). 
Actually, this holds true for $i=3$ as well:

\begin{proposition}[Hartl \cite{Hartl}]
\label{prop:low_degree}
For all N-series $N$ of $G$, one has that $D_3(N;\Z) = N_3$.
\end{proposition}

The most general result in the same direction is the following one:

\begin{theorem}[Mal'cev, Jennings, Hall, Lazard]
\label{th:MJHL}
Let $p$ be zero or a prime number. Let $N$ be a $p$-restricted N-series,
and assume that the ring $R$ has characteristic $p$.
Then, $D_i(N;R) = N_i$ for all $i\geq 1$.
\end{theorem}

\noindent
This statement is the compilation of two results.
In characteristic $p=0$, this is due to Mal'cev \cite{Mal'cev}, Hall \cite{Hall} and Jennings \cite{Jennings_bis} 
(when $N$ is the lower central series, at least) 
and, in characteristic $p>0$, this is due to Lazard \cite{Lazard}.
The reader is refered to Passi's book \cite[\S III]{Passi} for a detailed proof
solving, for an arbitrary N-series $N$, the two cases simultaneously.\\

\section{Homology cylinders and their finite-type invariants}

\label{sec:GH}

In this section, we start by recalling some generalities about the Goussarov--Habiro theory for $3$-manifolds
and, after that, we specialize to homology cylinders \cite{Goussarov,Habiro,GGP,Goussarov_bis}.
Next, we explain how the GHC can be regarded as an instance of the DSP.

\subsection{Finite-type invariants of $3$-manifolds}

\label{subsec:FTI}

First, we have to recall the notion of ``graph clasper.''
Start with a finite trivalent graph $\mathsf{G}$ decomposed as $\mathsf{G}_1 \cup \mathsf{G}_2$,
where $\mathsf{G}_1$ is a unitrivalent subgraph of $\mathsf{G}$ and 
$\mathsf{G}_2$ is a union of looped edges of $\mathsf{G}$. Give $\mathsf{G}$ a thickening\footnote{A \emph{thickening} of a graph 
$\mathsf{G}$ is a $\Z_2$-bundle over $\mathsf{G}$ with fiber $[-1,1]$.}
which is trivial on each looped edge of $\mathsf{G}_2$, and let $G$ be an embedding 
of this thickened graph into the interior of a compact oriented $3$-manifold $M$. 
Then, $G$ is said to be a \emph{graph clasper} in the manifold $M$.

The \emph{leaves} of $G$ are the framed knots in $M$ corresponding to the thickening of $\mathsf{G}_2$.
The \emph{shape} of $G$ is the abstract graph $\mathsf{G}_1$.
The \emph{degree} of $G$ is the number of trivalent vertices of its shape, which is assumed to be at least one.

\begin{example} If a graph clasper $G$ is $Y$-shaped, then it is called  a \emph{$Y$-graph}.
See Figure \ref{fig:Y} for the picture\footnote{By the ``blackboard framing'' convention, 
one-dimensional objects that are drawn on diagrams, 
such as graphs or links, are thickened along the plan.} of a $Y$-graph before embedding in a manifold $M$.
\begin{figure}[h!]
\begin{center}
\includegraphics[height=3cm,width=3cm]{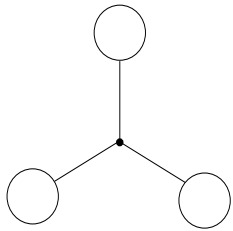} 
\caption{A $Y$-graph before embedding.}
\label{fig:Y}
\end{center}
\end{figure} 
\end{example}

\begin{figure}[h!]
\begin{center}
\includegraphics[height=4cm,width=4cm]{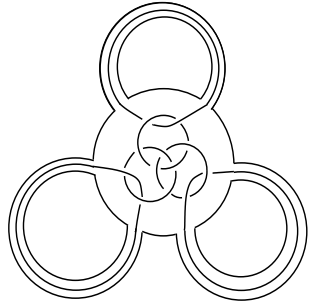} 
\caption{The framed link $B\subset \hbox{N}(G)$.}
\label{fig:B}
\end{center}
\end{figure}

A graph clasper carries surgery instructions to modify the manifold where it is embedded.
First, suppose that $G$ is a $Y$-graph in a manifold $M$ and let N$(G)$ be its
regular neighborhood in $M$. (This is a genus $3$ handlebody.)
The manifold obtained from $M$ by \emph{surgery along} $G$ is 
$$
M_G := \left(M\setminus \int\ \hbox{N}(G)\right) \cup \hbox{N}(G)_B
$$
where N$(G)_B$ is N$(G)$ surgered along the six-component framed link $B$ on Figure \ref{fig:B}.

\noindent
Let now $G$ be a graph clasper in $M$ of degree $k$. By applying the rule
\begin{center}
\includegraphics[height=1.5cm,width=8cm]{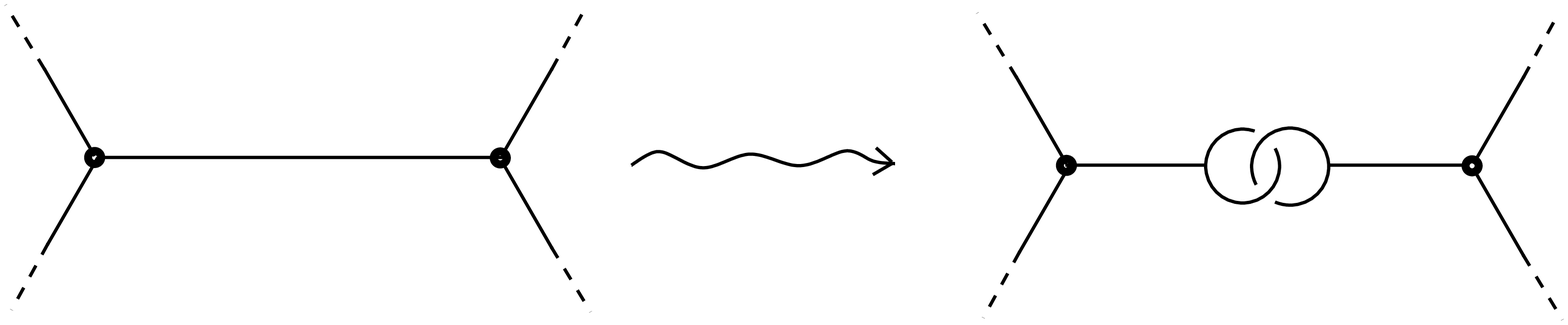} 
\end{center}
as many times as necessary, 
$G$ can be transformed to a disjoint union $Y(G)$ of $k$ $Y$-graphs in $M$.
The manifold obtained from $M$ by \emph{surgery along} $G$, denoted by $M_G$,
is the manifold $M$ surgered along each component of $Y(G)$.

\begin{definition}
\label{def:Y_k}
Let $k\in \N^*$. The \emph{$Y_k$-equivalence} is the equivalence relation 
among compact oriented $3$-manifolds generated 
by surgeries along connected graph claspers of degree $k$.
\end{definition}

We now fix a $Y_1$-equivalence class $\mathcal{M}_0$ and denote by $\Z\cdot \mathcal{M}_0$
the free Abelian group generated by its elements. 
Consider the following subgroup:
$$
\mathcal{F}_{k}^l\left(\mathcal{M}_0\right)
:= \left\langle\ [M,G]\ \left\vert\ M \in \mathcal{M}_0,\
G\subset M: \hbox{graph clasper with } |G|=l,
\deg(G)= k \right.\ \right\rangle_{\Z}
$$
where, regarding the graph clasper $G$ as the set of its connected components,
$|G|$ denotes its cardinality  and
$$
[M,G] := \sum_{G'\subset G} (-1)^{|G'|} \cdot M_{G'}\ \in \Z\cdot \mathcal{M}_0.
$$ 
One easily checks that
\begin{equation}
\label{eq:===}
\mathcal{F}_{k'}^l\left(\mathcal{M}_0\right)
\leq \mathcal{F}_{k}^l\left(\mathcal{M}_0\right) \leq
\mathcal{F}_{k}^{l'}\left(\mathcal{M}_0\right)
\quad \forall 1 \leq l \leq l' \leq k \leq k'.
\end{equation}
Then, one sets
$$
\mathcal{F}_d\left(\mathcal{M}_0\right)
:= \mathcal{F}_{d}^d\left(\mathcal{M}_0\right)
= \cup_{l=1}^d \mathcal{F}_{d}^l\left(\mathcal{M}_0\right)
= \cup_{k\geq d}\mathcal{F}_{k}^d\left(\mathcal{M}_0\right)
$$
to get the series
$$
\Z\cdot \mathcal{M}_0 \geq \mathcal{F}_1\left(\mathcal{M}_0\right)
\geq \mathcal{F}_2\left(\mathcal{M}_0\right) \geq \mathcal{F}_3\left(\mathcal{M}_0\right) \geq \cdots
$$
called the \emph{Goussarov--Habiro filtration}.

\begin{definition}
\label{def:fti}
Let $A$ be an Abelian group. A map $f:\mathcal{M}_0 \to A$ is a \emph{finite-type invariant} 
of \emph{degree} at most $d$ if its additive extension to $\Z\cdot \mathcal{M}_0$ 
vanishes on $\mathcal{F}_{d+1}\left(\mathcal{M}_0\right)$.
\end{definition}

The above definitions of a finite-type invariant and $Y_k$-equivalence 
are those which can be found\footnote{
The $Y_k$-equivalence is called $A_k$-equivalence in \cite{Habiro}, 
and $(k-1)$-equivalence in \cite{Goussarov}.} in Habiro's and Goussarov's papers. 
It is known \cite{Habiro} that those ``working'' definitions 
are equivalent to the ``expository'' definitions given in the introduction: 
This is proved in the Appendix.

As mentioned in the introduction, finite-type invariants of degree at most $d$ do not separate
two manifolds in the same $Y_{d+1}$-equivalence class:
This fact follows from the above definitions and (\ref{eq:===}). 
But the converse is not true for closed manifolds:

\begin{example}
\label{ex:counter-example}
Let us consider $M=\sharp^4\left( S^1 \times S^2\right)$ and 
$M'=(S^1\times S^1 \times S^1) \sharp (S^1 \times S^2)$. 
In contrast with the manifold $M'$, the manifold $M$
has a trivial triple-cup product form (say, with coefficients in $\Z$): 
Since the isomorphism class of the cohomology ring is preserved by surgery
along a connected degree $2$ graph clasper, the manifolds $M$ and $M'$ are not $Y_2$-equivalent.

Yet, they are not distinguished by degree $1$ finite-type invariants, which can be seen as follows.
The trivial $4$-component framed link $U$ in $S^3$ is a surgery presentation of $M$. 
Let $H$ and $G_1 \cup G_2$ be the following graph claspers in $M$ 
lying in $M\setminus U=S^3 \setminus U$:

\centerline{\relabelbox \small 
\epsfxsize 3.5truein \epsfbox{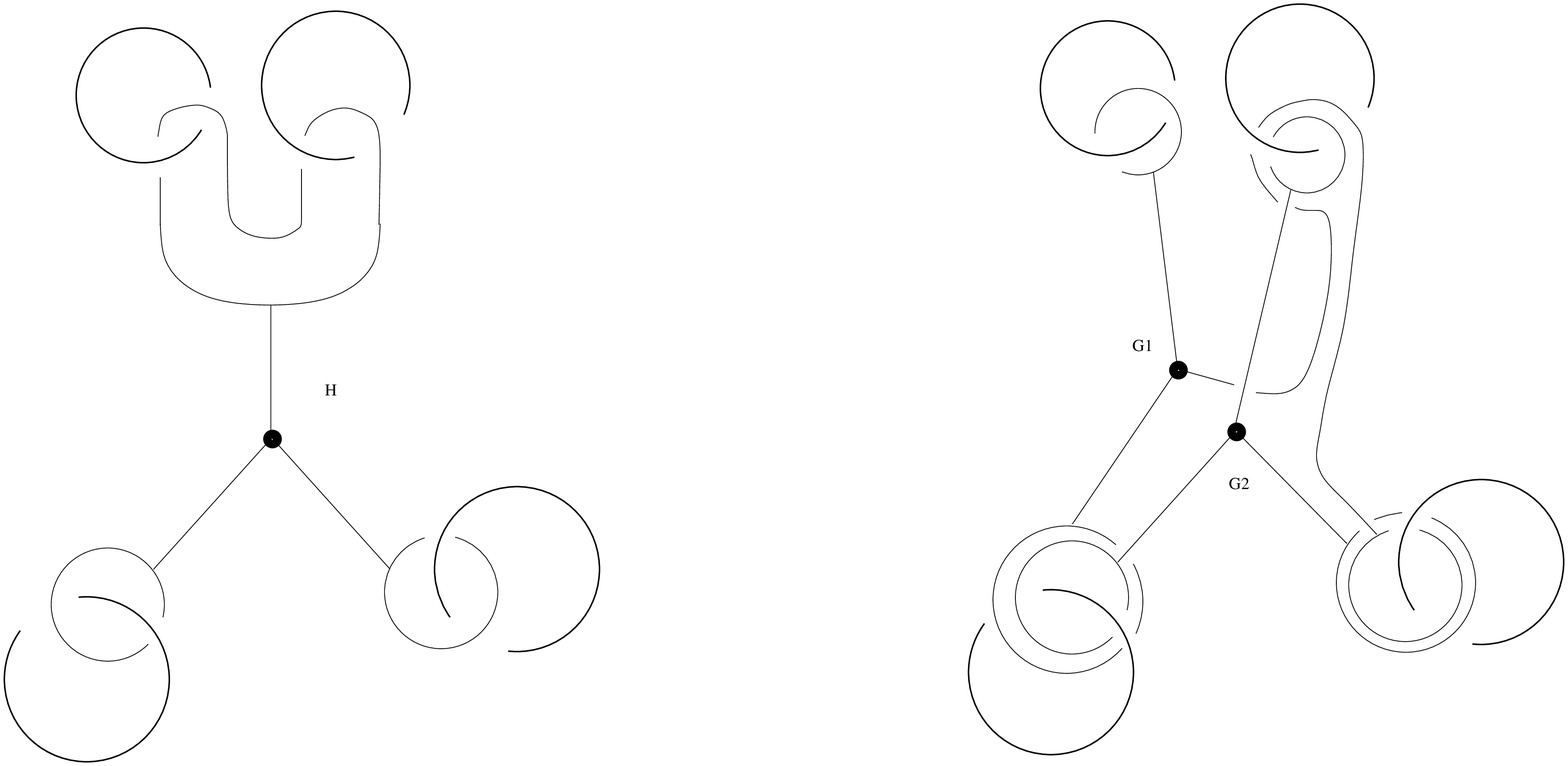}
\adjustrelabel <0.cm,-0.1cm> {H}{$H$}
\adjustrelabel <-0.3cm,-0.1cm> {G1}{$G_1$}
\adjustrelabel <-0.1cm,-0.2cm> {G2}{$G_2$}
\endrelabelbox}

\noindent
Observe that both $M_{G_1}$ and $M_{G_2}$ are homeomorphic to $M'$, since the latter can be presented
by surgery in $S^3$ along the Borromean rings union with a trivial $0$-framed knot.
Moreover, $M_{G_1\cup G_2}$ is seen to be homeomorphic to $M_H$ 
by the move 11 from \cite{Habiro} or, equivalently, by the move $Y_4$ from \cite{GGP}. 
But, $M_H$ is homeomorphic to $M'$ too,
as is clear if one of the two unknots grasped by the top leaf of $H$ is slided over its twin.
Hence, we conclude that the formal sum
$$
M - M'=  M- M_{G_1} - M_{G_2} + M_{G_1\cup G_2} = \left[M;G_1 \cup G_2 \right]
$$
belongs to $\mathcal{F}_2\left(\mathcal{M}_0\right)$,
where $\mathcal{M}_0$ is the $Y_1$-equivalence class of $M$ (or, equivalently, $M'$).
\end{example}

\subsection{Group of homology cylinders}

Let $\Sigma$ be a compact connected oriented surface, possibly with boundary. 

\begin{definition}
A \emph{cobordism} $M=(M,m)$ \emph{from} $\Sigma$ \emph{to} $\Sigma$ is a compact connected oriented $3$-manifold $M$
together with an orientation-preserving homeomorphism $m:\partial(\Sigma\times [-1,1]) \to \partial M$.
The \emph{trivial} cobordism from $\Sigma$ to $\Sigma$ is the cylinder $(\Sigma \times [-1,1],\Id)$.
A \emph{homology cylinder over} $\Sigma$ is a cobordism from $\Sigma$ to $\Sigma$ which is $Y_1$-equivalent to the trivial cobordism. 
\end{definition}

Let $\Cyl\left(\Sigma\right)$ be the set of homology cylinders,
up to homeomorphisms preserving the orientation and the boundary parametrization.
Given two homology cylinders $M$ and $M'$, we denote by $M'\cdot M$
the new one obtained by gluing the bottom of $M'$ to the top of $M$.
This ``stacking product''  turns $\Cyl\left(\Sigma\right)$ into a monoid,
whose identity element is the trivial cobordism $\Sigma\times [-1,1]$.

Homology cylinders can be recognized among cobordisms
thanks to the following homological characterization.

\begin{proposition}
\label{prop:Matveev}
A cobordism $(M,m)$ from $\Sigma$ to $\Sigma$ is a homology cylinder if, and only if, 
there exists an isomorphism $\psi:H_1(\Sigma \times [-1,1];\Z) \to H_1(M;\Z)$ 
such that the following diagram commutes:
$$
\xymatrix{
H_1\left(\partial(\Sigma \times [-1,1]);\Z\right) \ar[r]^{\incl_*} \ar[d]_{m_*}^\simeq& 
H_1\left(\Sigma \times [-1,1];\Z\right) \ar[d]^\psi_\simeq \\
H_1(\partial M;\Z) \ar[r]_{\incl_*} & H_1(M;\Z).
}
$$
\end{proposition}

\noindent
This fact has been announced by Habiro in \cite{Habiro}.
We propose below an alternative to the proof given by Habegger in \cite[\S 6]{Habegger}:

\begin{proof}
Assume that $(M,m)$ is a cobordism from $\Sigma$ to $\Sigma$, together with
an isomorphism $\psi : H_1(\Sigma \times [-1,1];\Z) \to H_1(M;\Z)$
such that the above diagram is commutative.
We wish to prove that $(M,m)$ is $Y_1$-equivalent to $(\Sigma \times [-1,1],\Id)$.
(The converse follows from the Mayer--Vietoris theorem,
since the $Y_1$-equivalence is generated by Torelli surgeries.)

First, assume that $\partial \Sigma$ is non-empty. Let $g$ be the genus of $\Sigma$
and let $r\geq 1$ be the number of its boundary components. 
Observe that $\Sigma\times [-1,1]$ is a handlebody of genus $n:= 2g+r-1$.
Let $x_1,\dots,x_n$ be closed oriented simple curves on $\partial \left(\Sigma \times [-1,1]\right)$
which are pairwise disjoint and whose images under $\incl_*$ generate $H_1\left(\Sigma \times [-1,1];\Z\right)$. 
Glue along each curve $x_i$ a $2$-handle: This produces a homology $3$-ball $B$ 
with, inside, an ordered oriented framed $n$-component tangle $\gamma$ 
whose components $\gamma_i$'s are the co-cores of the $2$-handles.
Let $D$ be a $3$-ball and let $\tau \subset D$ be an ordered oriented framed $n$-component tangle.
Identify $\partial B$ with $-\partial D$ in such a way that the endpoint of $\gamma_i$
goes to the startpoint of $\tau_i$, and vice-versa. Gluing $(B,\gamma)$ to $(D,\tau)$,
we obtain a homology $3$-sphere $\widehat{B}$ together with an ordered oriented framed $n$-component
link $\widehat{\gamma}$. Next, since $\partial M$ is identified with $\partial\left( \Sigma \times [-1,1]\right)$
via $m$, we can repeat \emph{exactly} the same constructions when $\Sigma \times [-1,1]$ is replaced by $M$: 
We get a compact oriented $3$-manifold $B'$ with a tangle $\gamma'$ and, next, 
a closed oriented $3$-manifold $\widehat{B'}$ with a link $\widehat{\gamma'}$.
The homological assumption on $m$ implies that $B'$ is a homology $3$-ball,
so that $\widehat{B'}$ is a homology $3$-sphere; it also implies that the linking matrices of 
$\widehat{\gamma'}$ and $\widehat{\gamma}$ are the same. 
Then, by combining \cite{Matveev} and \cite{MN}, we conclude that 
$\left(\widehat{B},\widehat{\gamma}\right)$ is $Y_1$-equivalent to 
$\left(\widehat{B'},\widehat{\gamma'}\right)$. 
This implies that $\left(B,\gamma\right)$ is $Y_1$-equivalent to 
$\left(B',\gamma'\right)$ or, which is the same, that 
$(\Sigma \times [-1,1],\Id)$ is $Y_1$-equivalent to $(M,m)$.

If now the surface $\Sigma$ is closed, we choose a disk $D\subset \Sigma$ and   
we consider the surface $\Sigma_0:= \Sigma \setminus\int\ D$. 
Let $T\subset M$ be a solid tube connecting $m(D\times (-1))$ to $m(D\times 1)$. 
We consider $M_0:= M \setminus \int\ T$ and we parametrize its boundary by a homeomorphism 
$m_0:\partial \left(\Sigma_0 \times [-1,1]\right) \to \partial M_0$ 
whose restriction to $\Sigma_0\times (-1) \cup \Sigma_0 \times 1$ coincides 
with the restriction of $m$. 
Using the fact that $H_1(\Sigma_0;\Z) \simeq H_1(\Sigma;\Z)$, one sees that
$\incl_*: H_1(\partial M_0;\Z) \to H_1(M_0;\Z)$ is surjective, and that
\begin{eqnarray*}
&&\Ker\left(\incl_*: H_1(\partial(\Sigma_0\times [-1,1]);\Z) 
\to H_1(\Sigma_0\times [-1,1];\Z) \right)\\
&=&\Ker \left(\incl_*\circ m_{0,*}: H_1(\partial(\Sigma_0\times [-1,1]);\Z) \to H_1(M_0;\Z)\right).
\end{eqnarray*}
So, by the previous paragraph, 
$(M_0,m_0)$ is $Y_1$-equivalent to $(\Sigma_0\times [-1,1],\Id)$. 
We conclude that $(M,m)$ is $Y_1$-equivalent to $(\Sigma\times [-1,1],\Id)$.
\end{proof}

Defining
$$
\Cyl_k\left(\Sigma\right)
:= \left\{M \in \Cyl\left(\Sigma\right) :
M \hbox{ is $Y_k$-equivalent to } \Sigma \times [-1,1] \right\},
$$
one gets a filtration of $\Cyl\left(\Sigma\right)$ by submonoids:
$$
\Cyl\left(\Sigma\right) = \Cyl_1\left(\Sigma\right) \geq \Cyl_2\left(\Sigma\right)
\geq \Cyl_3\left(\Sigma\right) \geq \cdots
$$
We recall some facts about that filtration:

\begin{theorem}[Goussarov, Habiro]
\label{th:GH}
For all $1 \leq k \leq l$, the quotient monoid
$\Cyl_k\left(\Sigma\right)/Y_l$ is a finitely generated group.
Moreover, for all $k,k'\geq 1$ and for all $l\geq k+k'$, the following inclusion holds:
$$
\left[\ \Cyl_{k}\left(\Sigma\right)/Y_l\ ,\
\Cyl_{k'}\left(\Sigma\right)/Y_l\ \right]\ \subset\
\Cyl_{k+k'}\left(\Sigma\right)/Y_l.
$$
\end{theorem}

\noindent
Those statements have been announced by Goussarov
in \cite[\S 4]{Goussarov} and by Habiro in \cite[\S 8.5]{Habiro}, but no proofs have been published.
Nevertheless, they proved the analogous results for string links
in full details in \cite[\S 9]{Goussarov_bis} and \cite[\S 5]{Habiro} respectively.
The reader may check that their proofs can be translated
to homology cylinders in a routine way.\footnote{ 
Indeed, some of their arguments get even simplified in the case of homology cylinders.}\\

The set $\Cyl(\Sigma)$ being a $Y_1$-equivalence class,
one can restrict the notion of finite-type invariants to homology cylinders 
(i.e$.$ take $\mathcal{M}_0:=\Cyl(\Sigma)$ in \S \ref{subsec:FTI}).
The next fact will play a key role in the sequel:

\begin{lemma}[Habiro]
\label{lem:Habiro}
For all integer $k\geq 1$, one has that
$$
\mathcal{F}_k\left(\Cyl(\Sigma)\right) =
\sum_{l=1}^k \sum_{\substack{k_1,\dots,k_l \geq 1\\ k_1 + \cdots + k_l =k}}
\mathcal{F}_{k_1}^1\left(\Cyl(\Sigma)\right)
\cdots \mathcal{F}_{k_l}^1\left(\Cyl(\Sigma)\right)
$$
in the monoid ring $\Z \cdot \Cyl(\Sigma)$.
\end{lemma}

\noindent
The inclusion ``$\supset$'' is easy to check. 
The converse inclusion ``$\subset$''  is proved by calculus of claspers
in a way parallel to the proof of \cite[Prop. 6.10]{Habiro},
which deals with string links.

\subsection{The Goussarov--Habiro conjecture}

Let $\Sigma$ be a compact connected oriented surface,
possibly with boundary.
We \emph{fix} an integer $k\geq 0$ and, in order to simplify notations, we set
$$
G := \Cyl\left(\Sigma\right)/Y_{k+1}
\quad \hbox{and} \quad N_i := \Cyl_i\left(\Sigma\right)/Y_{k+1}
\ \ \forall i\geq 1.
$$
By Theorem \ref{th:GH}, $N := (N_i)_{i\geq 1}$ is a finite N-series of the group $G$.

\begin{lemma} 
\label{lem:GH_to_DSP}
Let $d\geq 0$ be an integer and
let $R$ be a commutative ring. If two homology cylinders $M$ and $M'$ over $\Sigma$
are not distinguished by finite-type invariants of degree at most $d$
and with values in any $R$-module, then
$$
\{M\} = \{M'\} \ \mod \ D_{d+1}(N;R)
$$
where $\{M\}$ and $\{M'\}$ denote classes in $G$. 
Moreover, the converse is true if $d\leq k$.
\end{lemma}

\begin{proof}
The canonical projection $\Cyl\left(\Sigma\right) \to G$
induces a ring homomorphism
$$
R \cdot\Cyl\left(\Sigma\right) \to R[G]
$$
(from a monoid ring to a group ring).
By Lemma \ref{lem:Habiro},
this map sends $\mathcal{F}_{d+1}\left(\Cyl(\Sigma)\right) \otimes R$
to $\Delta_{d+1}(N;R)$,
and its kernel is $\mathcal{F}_{k+1}^1\left(\Cyl(\Sigma)\right)\otimes R$.
Hence the following commutative diagram:
$$
\xymatrix{
{\Cyl\left(\Sigma\right)} \ar[r] \ar[d]&
{R \cdot\Cyl\left(\Sigma\right) /
\mathcal{F}_{d+1}\left(\Cyl(\Sigma)\right)}\otimes R \ar[d]^-{\simeq\ \hbox{\scriptsize if } d\leq k}\\
{G}  \ar[r] & {R[G]/ \Delta_{d+1}(N;R)}.
}
$$
The above diagram shows that
\begin{eqnarray*}
M' - M \in  \mathcal{F}_{d+1}\left(\Cyl(\Sigma)\right) \otimes R
&\Longrightarrow &\{ M'\} - \{ M\} \in \Delta_{d+1}(N;R)\\
&\Longrightarrow& \{ M'\}\cdot \{ M\}^{-1} \in D_{d+1}(N;R)
\end{eqnarray*}
and that the converse implication holds true if $d\leq k$. 
Since the canonical map
$$
\Cyl\left(\Sigma\right) \to
{R \cdot\Cyl\left(\Sigma\right) /
\mathcal{F}_{d+1}\left(\Cyl(\Sigma)\right)}\otimes R
$$
is a degree $\leq d$ finite-type invariant with values in an $R$-module,
and since it is universal as such, the conclusions follow.
\end{proof}

Taking the coefficients in the ring $R=\Z$, we conclude that the GHC is a special instance of the DSP:

\begin{corollary}
Let $d\leq k$ be a non-negative integer.
The GHC is true in degree $d$ if, and only if, $D_{d+1}(N;\Z) = N_{d+1}$.
\end{corollary}

\begin{example}
Assume that $\Sigma$ is a $2$-disk. Then, $\Cyl(\Sigma)$ is the set of homology $3$-balls.
Since this monoid is Abelian, we deduce from Lemma \ref{lem:Abelian} 
that the GHC holds true for homology $3$-balls or, equivalently, for homology $3$-spheres.
The analogous result for knots has been proved 
by Goussarov and Habiro in \cite{Goussarov_bis,Habiro} where Lemma \ref{lem:Abelian} appears implicitely.
\end{example}

\begin{example}
We deduce from Proposition \ref{prop:low_degree} 
that the GHC is always true in degrees $d=0,1,2$.  
\end{example}

\noindent
That is the little we can say about integer coefficients.
Indeed, one can not apply Theorem \ref{th:MJHL} with coefficient ring $R=\Z$
since the above N-series $N$ is not $0$-restricted:
For instance, it is known that $\Cyl_1(\Sigma)/Y_2$ contains elements of order $2$
(detected by Birman--Craggs homomorphisms \cite{MM}).
Furthermore, for higher odd degree $i$, 
Habiro has found finite-order elements in $\Cyl_i(\Sigma)/Y_{i+1}$ \cite{Habiro_pc}.

\section{The Goussarov--Habiro conjecture with coefficients in a field}

\label{sec:weak}

As announced in the introduction, we now prove an analogue of the GHC 
for finite-type invariants with values in a fixed field 
(Theorem \ref{th:rational_GH} and Theorem \ref{th:modular_GH}).
Next, we deduce that the GHC holds true, at least, in a weaker form
(Corollary \ref{cor:big_degree}).

\subsection{Taking the coefficients in a field}

To prove Theorem \ref{th:rational_GH} and Theorem \ref{th:modular_GH},
we need the following two results respectively:

\begin{theorem}
\label{th:JH}
Let $G$ be a group, let $N$ be an N-series of $G$ and
let $\F$ be a field of characteristic zero.
Then, for all $i\geq 1$, one has that
$$
D_i(N;\F) = \sqrt{N_i}
$$
where $\sqrt{N_i} := \{g\in G : g^m \in N_i \hbox{ for some } m\geq 1\}$
is the isolator of $N_i$.
\end{theorem}

\begin{theorem}
\label{th:JL}
Let $G$ be a group, let $N$ be an N-series of $G$ and
let $\F$ be a field of characteristic $p>0$. 
Then, for all $i\geq 1$, one has that
$$
D_i(N;\F) = \prod_{\substack{k \geq 1,\ j\geq 0 \\ k\cdot p^j \geq i }} {N_k}^{(p^j)}
$$
where, for a subgroup $H\leq G$ and an integer $n\geq 1$, 
$H^{(n)}$ denotes the subgroup generated by all $n$-th powers of elements of $H$.
\end{theorem}

\noindent
Theorem \ref{th:JH} and Theorem \ref{th:JL} have been proved by Mal'cev, Jennings and Hall \cite{Mal'cev,Jennings_bis,Hall}
and by Jennings and Lazard \cite{Jennings,Lazard} respectively,
when $N$ is the lower central series of $G$.
But, the proofs given by Passi in his book \cite[\S IV]{Passi} 
can be adapted to any $N$-series of $G$. This is checked below.

Now, according to Lemma \ref{lem:GH_to_DSP}, 
Theorem \ref{th:modular_GH} follows immediately from Theorem \ref{th:JL}.
Similarly, Theorem \ref{th:rational_GH} is a consequence of Theorem \ref{th:JH} and the following lemma:

\begin{lemma}
Let $x$ and $y$ be elements of a nilpotent group $G$. 
There exists an $n\in \N^*$ such that $x^n=y^n$ if, and only if, 
$y^{-1} x$ has finite order.
\end{lemma}

\begin{proof}
Let $c$ be the nilpotency class of $G$, and assume that the lemma holds true for
groups of smaller nilpotency class.

Assume that $x=y\cdot z$, where $z\in G$ is such that $z^k=1$ for a certain $k>0$. 
Since the lemma holds true for the group $G/\gamma_{c}(G)$ (the length of its lower central series being $c-1$), 
one gets that $x^m =y^m \cdot u$ for some $m>0$ and $u\in \gamma_c(G)$.
Therefore, $u$ being central in $G$, $x^{mk} = y^{mk} \cdot u^k$. 
On the other hand, one has that $x^{mk} = y^{mk} \cdot z^{mk}$  modulo $[\Tors(G),G]$,
where $\Tors(G)$ denotes the subgroup of finite-order elements of $G$. Since $z^{mk}=1$, 
we deduce that $u^k \in [\Tors(G),G] \subset \Tors(G)$. 
So, there exists an $l>0$ such that $u^l=1$, from which we deduce that $x^{ml}=y^{ml}$.

Conversely, assume that $x^n = y^n$ for an $n>0$. Then, $\{x\}^n= \{y\}^n \in G/\Tors(G)$ as well: So,
if the converse were acquired for torsion-free nilpotent groups, 
we would have that $\{x\}= \{y\}\in G/\Tors(G)$ and we would be done. 
So, we can  assume that $G$ is torsion-free. 
By a theorem of Mal'cev (see \cite[\S 5.2.19]{Robinson} for instance), 
this implies that $G/\mathcal{Z}(G)$ is torsion-free as well,
where $\mathcal{Z}(G)$ denotes the center of $G$. Since the lemma holds true for $G/\mathcal{Z}(G)$
(the length of its upper central series being $c-1$), 
one obtains that $x=y \cdot z$ where $z$ is central. To the $n$-th power, 
this gives $x^n=y^n \cdot z^n$ from which we deduce that $z^n=1$ so that $z=1$ and $x=y$.  
\end{proof}

\subsubsection{Proof of Theorem \ref{th:JH}}

The proof given for the lower central series in \cite[\S IV.1.5]{Passi}
extends easily as follows:

\begin{lemma}
The sequence $\sqrt{N}:= \left(\sqrt{N_i}\right)_{i\geq 1}$
is a $0$-restricted $N$-series of $G$.
\end{lemma}

\begin{proof}
Since $N$ is an N-series,
$\gamma_i(G) \subset N_i$ so that $G/N_i$ is a nilpotent group.
So, the subset $\Tors(G/N_i)$ of finite-order elements of $G/N_i$ is a (normal) subgroup.
Since $\sqrt{N_i}$ is the inverse image of $\Tors(G/N_i)$ by the
canonical group homomorphism $G \to G/N_i$, we deduce that
$\sqrt{N_i}$ is a normal subgroup of $G$. 

Let $x\in \sqrt{N_m}$ and $y\in \sqrt{N_n}$ for some integers $m,n\geq 1$:
There exist $r,s\geq 1$ such that $x^r\in N_m$ and $y^s\in N_n$.
Let $\pi: G \to G/\sqrt{N_{n+m}}=:H$ be the canonical projection.
We are asked to show that $\pi([x,y])=1$. Assume the contrary.
Let $\left(\mathcal{Z}_i(H)\right)_{i\geq 0}$
denotes the upper central series of $H$.
There exists $i_0 \geq 1$ such that $\pi([x,y]) \in
\mathcal{Z}_{i_0}(H)\setminus \mathcal{Z}_{i_0-1}(H)$.
One easily checks, by recurrence on $k\geq 1$, that
$$
\forall \hbox{ group } K,\ \forall x,y \in K \hbox{ such that }[x,y]\in \mathcal{Z}(K), \
[x,y]^k = [x^k,y].
$$
So, $\pi([x,y])^{rs}=[\pi(x),\pi(y)]^{rs}$ is equal modulo $\mathcal{Z}_{i_0-1}(H)$
to $[\pi(x)^r,\pi(y)^s]=\pi([x^r,y^s])$.
Thus, modulo $\mathcal{Z}_{i_0-1}(H)$, $\pi([x,y])^{rs}$ belongs to
$\pi([N_m,N_n])\subset \pi(N_{m+n}) \subset \pi(\sqrt{N_{n+m}})=1$.
In other words, $\pi([x,y])$ gives a non-trivial
torsion element of the group $H/\mathcal{Z}_{i_0-1}(H)$.
But, by a theorem of Mal'cev \cite[\S 5.2.19]{Robinson},
$$
\forall \hbox{ nilpotent torsion-free group } K,\ \forall i\geq 0, \
K/\mathcal{Z}_i(K) \hbox{ is torsion-free}.
$$
The group $H$ is clearly torsion-free, and is nilpotent
since it is a quotient of $G/N_{n+m}$. So, we are led to a contradiction. 

Thus, $\sqrt{N}$ is an N-series, and it is clearly $0$-restricted. 
\end{proof}

\begin{lemma}
\label{lem:same}
For all integer $i\geq 1$, one has that 
$\Delta_i(\sqrt{N};\F)= \Delta_i(N;\F)$.
\end{lemma}

\begin{proof}
On the one hand, since $N_j \subset \sqrt{N_j}$ for all $j\geq 1$,
the inclusion $\Delta_i(N;\F) \subset \Delta_i(\sqrt{N};\F)$ certainly holds.

On the other hand, let $n\geq 1$ and $g\in \sqrt{N_n}$. There exists an $m\geq 1$ such that
$g^m \in N_n$, so that $g^m-1 \in \Delta_n(N;\F)$.
Besides, $(g-1)^n$ belongs to $I(\F[G])^n\subset \Delta_n(N;\F)$. 
But, it follows from the binomial identity (\ref{eq:binom_0}) that
$$
(g-1) \in (g^m-1) \cdot \F[G] + (g-1)^n\cdot \F[G].
$$
So, $(g-1)$ belongs to $\Delta_n(N;\F)$.
We conclude that $\Delta_i(\sqrt{N};\F) \subset \Delta_i(N;\F)$.
\end{proof}

Thus, Theorem \ref{th:JH} follows from Theorem \ref{th:MJHL} 
applied to the $0$-restricted N-series $\sqrt{N}$.

\subsubsection{Proof of Theorem \ref{th:JL}}

Again, the proof given for the lower central series 
in \cite[\S IV.1.9]{Passi} extends directly as follows: Consider
\begin{equation}
\label{eq:Lazard}
L_i(N;p):= \prod_{\substack{k \geq 1,\ j\geq 0 \\ k\cdot p^j \geq i }} {N_k}^{(p^j)}.
\end{equation}
This is a normal subgroup of $G$.
The sequence $L(N;p):=\left\{L_i(N;p)\right\}_{i\geq 1}$ is called the \emph{Lazard series}
defined by the N-series $N$ in characteristic $p$.

\begin{lemma}
The sequence $L(N;p)$ is a $p$-restricted N-series of $G$.
\end{lemma}

\begin{proof}
The proof is exactly the same as for \cite[\S IV.1.22]{Passi}, 
with each occurrence of $\gamma_i(G)$ replaced by $N_i$.
\end{proof}

\begin{lemma}
For all integer $i\geq 1$, one has that $\Delta_i(L(N;p);\F)=\Delta_i(N;\F)$.
\end{lemma}

\begin{proof}
Since $N_j \subset L_j(N;p)$ for all $j\geq 1$, we have that 
$\Delta_i(N;\F) \subset \Delta_i(L(N;p);\F)$. 

For all $j\geq 1$, for all $k\geq 1$ and $l\geq 0$ such that $kp^l\geq j$ and for all $x\in N_k$,
one has that $x\in D_k(N;\F)$ (since $N_k \subset D_k(N;\F)$) so that 
$x^{p^l} \in D_{kp^l}(N;\F) \subset D_j(N;\F)$ 
(since $D(N;\F)$ is a $p$-restricted N-series by Lemma \ref{lem:filtrations}). 
It follows that $L_j(N;p)\subset D_j(N;\F)$ for all $j\geq 1$. We conclude that
$\Delta_i(L(N;p);\F) \subset \Delta_i(D(N;\F);\F)=\Delta_i(N;\F)$.  
\end{proof}

So, Theorem \ref{th:JL} follows from Theorem \ref{th:MJHL} 
applied to the $p$-restricted N-series $L(N;p)$. 

\subsection{Increasing the degree}

In view of Lemma \ref{lem:GH_to_DSP}, Corollary \ref{cor:big_degree} 
is an application of the following statement:

\begin{corollary}
\label{cor:consequence_of_JH_and_JL}
Let $G$ be a finitely generated group and let $N$ be an $N$-series of $G$ such that 
$N_{d+1}=\{1\}$ for some $d\geq 0$. Set
$$
m := \max\left\{ \left. p^{e(p)} \right\vert p \in \mathcal{P} \right\}
$$
where $\mathcal{P}$ denotes the set of prime numbers $p$ 
and $e(p)$ is the $p$-exponent of the group $G$.
Then, $m$ is finite and 
the dimension subgroup $D_{m^d(d+1)}(N;\Z)$ is trivial.
\end{corollary}

\noindent
This is proved using Theorem \ref{th:JH} and Theorem \ref{th:JL} as follows:

\begin{proof}[Proof of Corollary \ref{cor:consequence_of_JH_and_JL}]
Since the N-series $N$ collapses at the $(d+1)$-st level, $G$ is a nilpotent group of class at most $d$.
So, the set $\Tors(G)$ of its finite-order elements is a subgroup 
and is a direct product indexed by $p\in\mathcal{P}$ of $p$-groups:
$$
\Tors(G) = \bigodot_{p \in \mathcal{P}} \Tors_p(G) \ \hbox{ where } \ 
\Tors_p(G) := \left\{x\in G: x^{p^e}=1 \hbox{ for some } e\geq 0 \right\}.
$$
(See, for instance, \cite[\S 5.2.7]{Robinson}.) 
Recall that the $p$-exponent $e(p)$ of $G$ is then defined to be the smallest $e$ 
such that $x^{p^e}=1$ for all $x\in \Tors_p(G)$.
Next, $G$ being finitely generated and nilpotent, $\Tors(G)$ is finite. 
So, for each $p \in \mathcal{P}$, $e(p)$ is finite and $e(p)$ is non-zero
only for finitely many $p$. We conclude that the number $m$ is finite.

Let $x \in D_{m^d(d+1)}(N;\Z)$. We wish to show that $x=1$.
Since $m^d(d+1) \geq d+1$, $x$ belongs to $D_{d+1}(N;\Q)$ 
which coincides with $\Tors(G)$ by Theorem \ref{th:JH}. 
Therefore, $x$ can be written uniquely as a finite product  
$$
x=\prod_{p\in \mathcal{P}} x_p \quad \hbox{ where } x_p \in \Tors_p(G).
$$
Let us fix a prime number $p$ and let us prove that the coordinate $x_p$ of $x$ is trivial.

Since $m^d(d+1) \geq p^{d\cdot e(p)}(d+1)$, $x$ belongs to $D_{p^{d\cdot e(p)}(d+1)}(N;\Z/p\Z)$ 
which coincides with
$$
\prod_{\substack{k \geq 1,\ j\geq 0 \\ k\cdot p^j \geq  (d+1) \cdot p^{d\cdot e(p)}}} {N_k}^{(p^j)}
$$
according to Theorem \ref{th:JL}. Then, using the hypothesis that $N_{d+1}=1$, one can write 
$$
x= y_1^{p^{d \cdot e(p)}} \cdots y_r^{p^{d \cdot e(p)}}
\quad \hbox{where} \ y_1,\dots,y_r \in G.
$$
By a result of Mal'cev (see \cite[Lemma 2.4.2]{Hartley}), 
for all $r\geq 1$ and in every nilpotent group of class $d$,
every product of $r^d$-th powers  is an $r$-th power. So, one can write
$$
x=y^{p^{e(p)}} \hbox{ where } y\in G.
$$ 
Since $x$ has finite order, $y$ has too. Thus, $y=y_p \cdot Y $ 
where $y_p\in \Tors_p(G)$ and $Y$ belongs to $\odot_{q\neq p} \Tors_q(G)$. 
We conclude that
$$
x= {y_p}^{p^{e(p)}} \cdot Y^{p^{e(p)}} = Y^{p^{e(p)}} \quad \in \bigodot_{q\neq p} \Tors_q(G)
$$
so that $x_p$ must be trivial.
\end{proof}

\section{The algebra dual to finite-type invariants of homology cylinders}

\label{sec:Quillen}

We start by extending a result due to Quillen \cite{Quillen}
for the lower central series of a group, to arbitrary N-series.
Applications to homology cylinders follow.

\subsection{Quillen's theorem for an arbitrary N-series}

Let $G$ be a group and let $N$ be an N-series of $G$. 
The following construction is due to Lazard \cite{Lazard}:

\begin{definition}
The graded Lie ring \emph{induced} by $N$ is the graded Abelian group
$$
\mathcal{L}(N):= \bigoplus_{i\geq 1} \frac{N_i}{N_{i+1}}
$$
with the Lie bracket\footnote{This is well-defined and satisfies the axioms of a Lie bracket
thanks to the Hall--Witt identities.} defined by
$$
\left[\{n_i\},\{n_j\}\right] := \{\left[n_i,n_j\right]\} \quad \in N_{i+j}/N_{i+j+1}
$$
on homogeneous elements $\{n_i\} \in N_{i}/N_{i+1}$ and $\{n_j\} \in N_{j}/N_{j+1}$.
\end{definition}

Now, let $R$ be a commutative ring.

\begin{definition}
The graded $R$-algebra \emph{induced} by $N$ is the graded Abelian group
$$
\mathcal{A}(N;R) := \bigoplus_{i\geq 0} \frac{\Delta_i(N;R)}{\Delta_{i+1}(N;R)}
\quad \hbox{where } \Delta_0(N;R):=R[G]
$$
with the multiplication coming from the ring $R[G]$.
\end{definition}

There is a graded group homomorphism
$$
\theta: \mathcal{L}(N) \to \mathcal{A}(N;R)
$$
defined by $\{n_i\}  \mapsto \{n_i-1\}$ 
for all homogeneous element $\{n_i\} \in N_i/N_{i+1}$.

If $\mathcal{A}(N;R)$ is given the Lie bracket defined by its associative multiplication,
then $\theta\otimes R$ is a morphism of Lie $R$-algebras 
since, for all $\{n_i\} \in N_i/N_{i+1}$ and $\{n_j\} \in N_j/N_{j+1}$,
$$
(\theta\otimes R)\left(\left[\{n_i\}\otimes 1,\{n_j\}\otimes 1\right]\right) 
= (\theta\otimes R)(\{[n_i,n_j]\}\otimes 1)
$$
$$
\stackrel{\footnotesize{(\ref{eq:comm})}}{=} \left\{[n_i-1,n_j-1] + (n_i^{-1}n_j^{-1}-1)[n_i-1,n_j-1]\right\}
$$
$$
= \left\{[n_i-1,n_j-1] \right\} = \left[\{n_i-1\},\{n_j-1\}\right] = 
\left[(\theta\otimes R)(\{n_i\}\otimes 1),(\theta\otimes R)(\{n_j\}\otimes 1)\right].
$$
Thus, passing to the universal enveloping algebra, one gets an $R$-algebra homomorphism
$$
\U(\theta \otimes R) : \U\left(\mathcal{L}(N) \otimes R\right) \to \mathcal{A}(N;R)
$$
which is clearly surjective.\\

We now specialize to the case when $R$ is a field $\F$:

\begin{theorem}
\label{th:Quillen}
Let $\F$ be a field of characteristic $p\geq 0$,
and let $N$ be a $p$-restricted N-series of a group $G$. 
If $p=0$, then the $\F$-algebra homomorphism
$$
\U(\theta \otimes \F): \U\left(\mathcal{L}(N)\otimes \F\right) \to \mathcal{A}(N;\F)
$$
is an isomorphism. If $p>0$, then $\U(\theta \otimes \F)$ factors to an $\F$-algebra isomorphism
$$
\V(\theta \otimes \F): \V\left(\mathcal{L}(N)\otimes \F\right) \to \mathcal{A}(N;\F)
$$
where $\V(-)$ denotes the \emph{restricted} enveloping algebra functor.
\end{theorem}

\noindent
When $N=\gamma(G)$, this is Quillen's theorem \cite{Quillen}. 
The same strategy of proof works for any N-series $N$: 

\begin{proof}
First of all, it follows from Theorem \ref{th:MJHL} that
$\theta:  \mathcal{L}(N) \to \mathcal{A}(N;\F)$ is injective,
in both cases $p>0$ and $p=0$. This implies that 
$\theta \otimes \F_p: \mathcal{L}(N)\otimes \F_p \to \mathcal{A}(N;\F_p)$ is injective as well,
where $\F_p$ denotes $\Z/p\Z$ if $p>0$ and $\Q$ if $p=0$. 
So, $\left(\theta \otimes \F_p\right) \otimes_{\F_p} \F$ is injective too. 
There is the commutative diagram
$$
\xymatrix{
{\left(\L(N) \otimes \F_p\right) \otimes_{\F_p} \F}  \ar[r]^-\simeq 
\ar[d]_-{(\theta \otimes \F_p)\otimes_{\F_p} \F} & 
{\L(N) \otimes \F} \ar[d]^{\theta \otimes\F}\\
{\A(N;\F_p)\otimes_{\F_p} \F} \ar[r] & {\A(N;\F),}
}
$$
whose bottom map is bijective since
$$
\frac{\Delta_i(N;\F_p)}{\Delta_{i+1}(N;\F_p)} \otimes_{\F_p}\F \simeq
\frac{\Delta_i(N;\Z)\otimes \F_p}{\Delta_{i+1}(N;\Z)\otimes \F_p} \otimes_{\F_p}\F \simeq
\frac{\left(\Delta_i(N;\Z)\otimes \F_p\right)\otimes_{\F_p}\F}
{\left(\Delta_{i+1}(N;\Z)\otimes \F_p\right)\otimes_{\F_p}\F}
$$ 
$$
\simeq
\frac{\Delta_i(N;\Z)\otimes \F}{\Delta_{i+1}(N;\Z)\otimes \F} \simeq 
\frac{\Delta_i(N;\F)}{\Delta_{i+1}(N;\F)}.
$$
We conclude that, in general,
\begin{equation}
\label{eq:injectivity}
\hbox{$\theta\otimes \F:  \mathcal{L}(N)\otimes \F \to \mathcal{A}(N;\F)$ is injective.}
\end{equation}

In the case when $p>0$, we define a map
$$
\xi : \mathcal{L}(N)\otimes \F \to \mathcal{L}(N) \otimes \F
$$ 
by $\xi(\{n_i\} \otimes 1) = \{n_i^p\} \otimes 1 \in N_{ip}/N_{ip+1}$ 
on a homogeneous element $\{n_i\} \in N_i/N_{i+1}$. Observe that
$$
(\theta\otimes \F)\left(\{n_i\} \otimes 1\right)^{\ p} = \{n_i-1\}^p = 
\{(n_i-1)^p\} \stackrel{\footnotesize{(\ref{eq:binom_p})}}{=} 
\{n_i^p-1\}  = (\theta\otimes \F)\circ \xi(\{n_i\}).
$$
This observation together with (\ref{eq:injectivity})
imply, first, that $\mathcal{L}(N)\otimes \F$ endowed with $\xi$ is a restricted Lie $\F$-algebra
and, second, that $\theta\otimes \F$ is a restricted Lie $\F$-algebra morphism. 
Therefore, in the case $p>0$, the map $\U(\theta \otimes \F)$ factors to an $\F$-algebra epimorphism
$\V(\theta \otimes \F): \V\left(\mathcal{L}(N)\otimes \F\right) \to \mathcal{A}(N;\F)$.

Let $\delta: \F[G] \to \F[G] \otimes_\F \F[G]$ be the usual coproduct of $\F[G]$, 
defined by $\delta(g)= g \otimes g$ for all $g\in G$. Since
\begin{equation}
\label{eq:coproduct}
\delta(x-1)= 1 \otimes (x-1) + (x-1) \otimes 1 + (x-1) \otimes (x-1)
\end{equation}
for all $x\in G$, one obtains that
$$
\delta\left(\Delta_n(N;\F)\right) \subset 
\sum_{i+j=n} \Delta_i(N;\F) \otimes_\F \Delta_j(N;\F).
$$
Thus, the (co-commutative) Hopf algebra structure of $\F[G]$ induces 
a (co-commutative) Hopf algebra structure on $\mathcal{A}(N;\F)$. 
It follows also from (\ref{eq:coproduct}) that, for all $i\geq 1$ and for all $n_i \in N_i$,
the homogeneous element $\{n_i-1\}$ is primitive in $\mathcal{A}(N;\F)$. 
Thus, $\mathcal{A}(N;\F)$ is generated by its set of primitives $\P\mathcal{A}(N;\F)$ 
and, by Milnor--Moore theorem \cite{MM},
\begin{equation}
\label{eq:MM}
\mathcal{A}(N;\F) = \left\{ 
\begin{array}{l}
\U \P \mathcal{A}(N;\F) \quad \hbox{if } p=0\\
\V \P \mathcal{A}(N;\F) \quad \hbox{if } p>0.\\
\end{array}\right.
\end{equation}
Since $\U(\theta\otimes \F)$ and $\V(\theta\otimes \F)$ are surjective for $p=0$ and $p>0$ respectively, 
the Poincar\'e--Birkhoff--Witt theorem implies that $\theta\otimes \F$ is surjective onto $\P\mathcal{A}(N;\F)$ in both cases. 
We then deduce from (\ref{eq:injectivity}) that 
$\theta\otimes \F: \mathcal{L}(N) \otimes \F \to \P \mathcal{A}(N;\F) $ is an isomorphism, 
and we conclude by (\ref{eq:MM}).
\end{proof}

\begin{corollary}
\label{cor:U}
Let $\F$ be a field of characteristic zero and let $N$ be an N-series of a group $G$.
Then, the $\F$-algebra morphism
$$
\U(\theta\otimes \F): \U\left(\mathcal{L}(N)\otimes \F\right) \to \mathcal{A}(N;\F)
$$
is an isomorphism.
\end{corollary}

\begin{proof}
It is enough to apply Theorem \ref{th:Quillen} to the $0$-restricted series $D(N;\F)$ which,
by Theorem \ref{th:JH}, coincides with $\sqrt{N}$, and to observe that the canonical map
$N_i/N_{i+1} \otimes \F \to \sqrt{N_i}/\sqrt{N_{i+1}} \otimes \F$
is an isomorphism. 
\end{proof}

\begin{corollary}
\label{cor:V}
Let $\F$ be a field of characteristic $p>0$ and let $N$ be an N-series of a group $G$.
Then, the $\F$-algebra morphism
$$
\V(\theta\otimes \F): \V\left(\mathcal{L}(L(N;p))\otimes \F\right) \to \mathcal{A}(N;\F),
$$
where $L(N;p)$ denotes the Lazard series defined by (\ref{eq:Lazard}), is an isomorphism.
\end{corollary}

\begin{proof}
It is enough to apply Theorem \ref{th:Quillen} to the $p$-restricted series $D(N;\F)$ which,
by Theorem \ref{th:JL}, coincides with $L(N;p)$. 
\end{proof}

\subsection{Application to homology cylinders}

Let $\Sigma$ be a compact connected oriented surface.
Applying the above construction to the group $G:= \Cyl(\Sigma)/Y_{k+1}$ 
and to its N-series $N:=\left(\Cyl_i(\Sigma)/Y_{k+1}\right)_i$ for all $k\geq 1$, 
one gets the \emph{graded Lie ring of homology cylinders} over $\Sigma$:
$$
\overline{\Cyl}(\Sigma):= \bigoplus_{i\geq 1} \frac{\Cyl_i(\Sigma)}{Y_{i+1}}.
$$
It has been introduced by Habiro in \cite[\S 8.5]{Habiro}. 

\begin{corollary}
Let $\F$ be a field of characteristic zero.
There is an $\F$-algebra isomorphism
$$
\U\left(\overline{\Cyl}(\Sigma) \otimes \F \right) 
\simeq \bigoplus_{i\geq 0} 
\frac{\mathcal{F}_i\left(\Cyl(\Sigma)\right)}{\mathcal{F}_{i+1}\left(\Cyl(\Sigma)\right)} \otimes \F
$$
defined by $\{M\}\otimes 1 \mapsto \{M-\Sigma\times [-1,1]\}$ for all $M \in \Cyl_i(\Sigma)$.
\end{corollary}

Let $p>0$ be a prime number. According to the previous subsection, the graded Abelian group
$$
\overline{\Cyl}(\Sigma;p) := 
\bigoplus_{i\geq 1} \frac{\left. \left(\prod_{kp^j\geq i}\Cyl_k(\Sigma)^{(p^j)}\right) \right/ Y_{i+1}}
{\left. \left(\prod_{kp^j\geq i+1}\Cyl_k(\Sigma)^{(p^j)}\right) \right/ Y_{i+1}}
$$
is a restricted Lie $\Z/p\Z$-algebra.

\begin{corollary}
Let $\F$ be a field of characteristic $p>0$.
There is an $\F$-algebra isomorphism
$$
\V \left(\overline{\Cyl}(\Sigma;p) \otimes \F\right)
\simeq \bigoplus_{i\geq 0} 
\frac{\mathcal{F}_i\left(\Cyl(\Sigma)\right)}{\mathcal{F}_{i+1}\left(\Cyl(\Sigma)\right)} \otimes \F
$$
defined by $\{M\}\otimes 1 \mapsto \{M-\Sigma\times [-1,1]\}$ 
for all $M \in \prod_{kp^j\geq i}\Cyl_k(\Sigma)^{(p^j)}$.
\end{corollary}

\section{Appendix: Torelli groups and claspers}

The following two lemmas explain why the definitions of a finite-type invariant and of the $Y_k$-equivalence
given in \S \ref{subsec:FTI} agree with those exposed in the introduction. 

\begin{lemma}
\label{lem:equivalence_FTI}
Let $\mathcal{M}_0$ be a $Y_1$-equivalence class of compact oriented $3$-manifolds. 
The $(d+1)$-st subgroup $\mathcal{F}_{d+1}(\mathcal{M}_0)$
of the Goussarov--Habiro filtration (as defined at \S \ref{subsec:FTI}) 
is generated by alternate sums
$$
\sum_{\Gamma' \subset \Gamma} (-1)^{|\Gamma'|}\cdot M_{\Gamma'} \quad \in \Z\cdot \mathcal{M}_0
$$
where $M\in \mathcal{M}_0$, $\Gamma$ is a set of $(d+1)$ pairwise disjoint handlebodies in $M$
with, specified for each, a Torelli automorphism of the boundary,
and where $M_{\Gamma'}$ is obtained from $M$ by Torelli surgery along those handlebodies belonging to $\Gamma'$.
\end{lemma}

\noindent 
This is essentially \cite[Theorem 5.2]{GGP}.

\begin{proof}[Proof of Lemma \ref{lem:equivalence_FTI}]
According to \cite{Johnson}, the Torelli group of a closed oriented surface $S$ (of genus at least $3$) 
is generated by \emph{bounding-pair maps}, i.e$.$ simultaneous opposite Dehn twists $T_\alpha^{-1} \circ T_\beta$ where
$(\alpha,\beta)$ is a pair of simple closed curves that bound a genus one sub-surface $\Sigma$ of $S$.
See Figure \ref{fig:genus_one}. 
Moreover, the mapping cylinder of $T_\alpha^{-1} \circ T_\beta$, 
seen as a cobordism from $\Sigma$ to $\Sigma$, is obtained from the trivial cobordism $\Sigma \times [-1,1]$
by surgery along a $Y$-graph: 
This follows from Lickorish's trick -- which converts Dehn twists into surgery along knots -- 
and calculus of claspers as shown on Figure \ref{fig:calculus}.

This proves that the surgery along a $Y$-graph can be realized by a Torelli surgery and that, conversely,
any Torelli surgery can be decomposed into a sequence of surgeries along $Y$-graphs. 
Since $\mathcal{F}_{d+1}(\mathcal{M}_0) = \mathcal{F}_{d+1}^{d+1}(\mathcal{M}_0)$, 
the conclusion easily follows.

\begin{figure}[h!]
\centerline{\relabelbox \small 
\epsfxsize 2truein \epsfbox{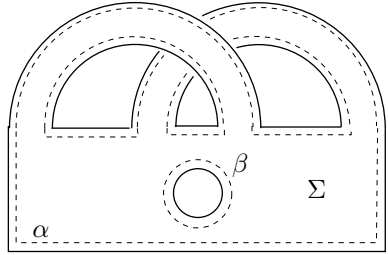}
\adjustrelabel <0.cm,-0.1cm> {Sigma}{$\Sigma$}
\adjustrelabel <-0.cm,-0.cm> {a}{$\alpha$}
\adjustrelabel <-0.cm,-0.cm> {b}{$\beta$}
\endrelabelbox}
\caption{A genus one surface $\Sigma$ with two boundary components $\alpha$ and $\beta$.}
\label{fig:genus_one}
\end{figure}
\begin{figure}[h!]
\centerline{\relabelbox \small 
\epsfxsize 5.5truein \epsfbox{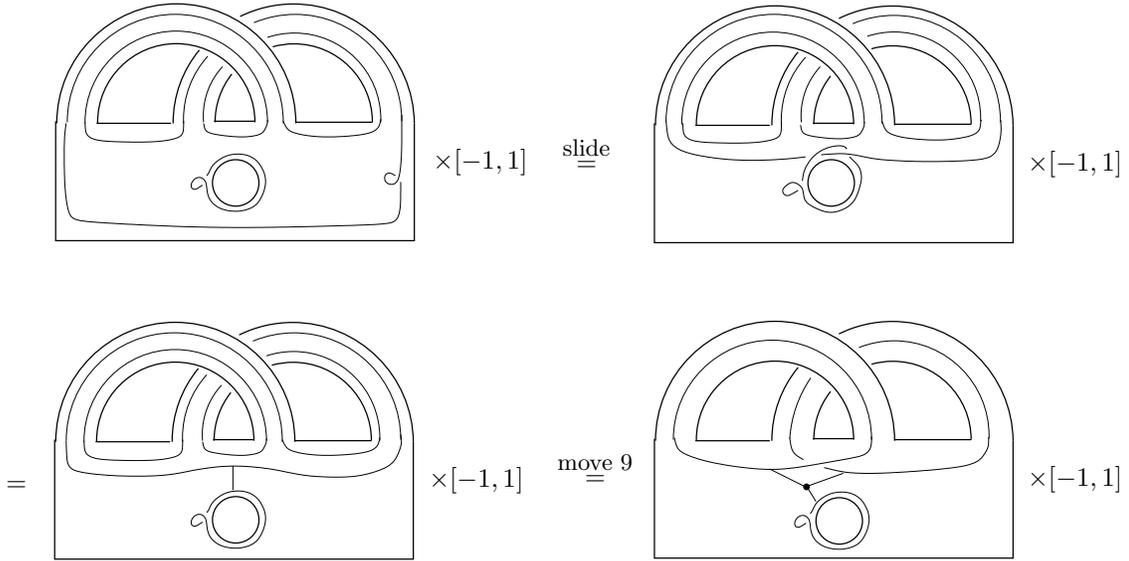}
\adjustrelabel <0.cm,-0.1cm> {times1}{$\times [-1,1]$}
\adjustrelabel <0.cm,-0.1cm> {times2}{$\times [-1,1]$}
\adjustrelabel <0.cm,-0.1cm> {times3}{$\times [-1,1]$}
\adjustrelabel <0.cm,-0.1cm> {times4}{$\times [-1,1]$}
\adjustrelabel <0.4cm,-0.1cm> {=1}{$\stackrel{\footnotesize{\hbox{slide}}}{=}$}
\adjustrelabel <-0.2cm,-0.1cm> {=2}{$=$}
\adjustrelabel <0.3cm,-0.1cm> {=3}{$\stackrel{\footnotesize{\hbox{move 9}}}{=}$}
\endrelabelbox}
\caption{Equivalent surgery descriptions of the mapping cylinder of a bounding-pair map
(using calculus of claspers \cite{Habiro}).}
\label{fig:calculus}
\end{figure}
\end{proof}

\begin{lemma}
\label{lem:equivalence_Yk}
Let $M$ and $M'$ be compact oriented $3$-manifolds and let $k\in  \N^*$.
Then, $M$ and $M'$ are $Y_k$-equivalent (in the sense of \S \ref{subsec:FTI}) if, and only if,
there exists a handlebody $H\subset M$ and an element $h$
of the $k$-th term of the lower central series of the Torelli group of $\partial H$
such that the Torelli surgery $M\leadsto M_h$ produces $M'$.
\end{lemma}

\noindent
This has been announced by Habiro in \cite{Habiro}.
We are indebted to him for a great simplification in the proof below.

\begin{proof}[Proof of Lemma \ref{lem:equivalence_Yk}]
Let $\mathcal{I}(\Sigma)$ denote the Torelli group of a closed oriented surface $\Sigma$.
We deduce from Proposition \ref{prop:Matveev} that the ``mapping cylinder'' map
$$
h \mapsto \left(\Sigma \times [-1,1], (h\times (-1)) \cup (\Id \times 1) \right) 
$$
transforms a Torelli automorphism of $\Sigma$ to a homology cylinder over $\Sigma$, hence a monoid homomorphism
$m: \mathcal{I}(\Sigma) \to \Cyl(\Sigma)$. Let $k\geq 1$ be an integer.
By composing $m$ with the canonical projection $\Cyl(\Sigma) \to \Cyl(\Sigma)/Y_{k}$,
one gets a group homomorphism which, by Theorem \ref{th:GH},
sends $\gamma_k(\mathcal{I}(\Sigma))$ to $\gamma_k(\Cyl(\Sigma)/Y_k)=1$. 
Thus, $m(\gamma_k(\mathcal{I}(\Sigma)))$ is contained in $\Cyl_k(\Sigma)$,
and this proves the implication ``$\Leftarrow$''.

Now, let us show the converse implication ``$\Rightarrow$''. Let $k\geq 1$ be an integer.
On the one hand, the relation between two compact oriented $3$-manifolds $M$ and $M'$
\begin{quote}
`` There exists a handlebody $H\subset M$ and an  $h \in \gamma_k(\mathcal{I}(\partial H))$
such that the Torelli surgery $M\leadsto M_h$ gives $M'$. ''
\end{quote}
is transitive\footnote{Indeed, any two handlebodies embedded in a $3$-manifold can be isotoped to disjoint positions 
and, next, they can be connected by a solid tube.} and is symmetric.  
On the other hand, it is easily seen from the definition that the $Y_k$-equivalence is generated
by surgeries along \emph{tree} claspers, i.e$.$ graph claspers whose shape is a tree. 
So, it is enough to prove that, given a compact oriented $3$-manifold $M$ 
and given a connected tree clasper $G \subset M$ of degree $k$, 
$M_G$ can be obtained from $M$ by a Torelli surgery whose
gluing homeomorphism lives in the $k$-th term of the lower central series. 
Let $\hbox{N}(G)$ be a regular neighborhood of $G$ in $M$. This is a genus $(k+2)$ handlebody
which we identify with $[-1,1]^3 \setminus\int\ \hbox{N}(\tau)$, where 
$\tau$ is the $(k+2)$-component framed trivial tangle of the cube drawn on
the left-hand part of Figure \ref{fig:picture}. 
According to \cite[Lemma 3.20]{Habiro}, surgery along $G \subset [-1,1]^3$
modifies $\tau$ by insertion of a framed pure braid $\beta$, 
as shown on the right-hand part of the same figure.
\begin{figure}
\centerline{\relabelbox \small 
\epsfxsize 6truein \epsfbox{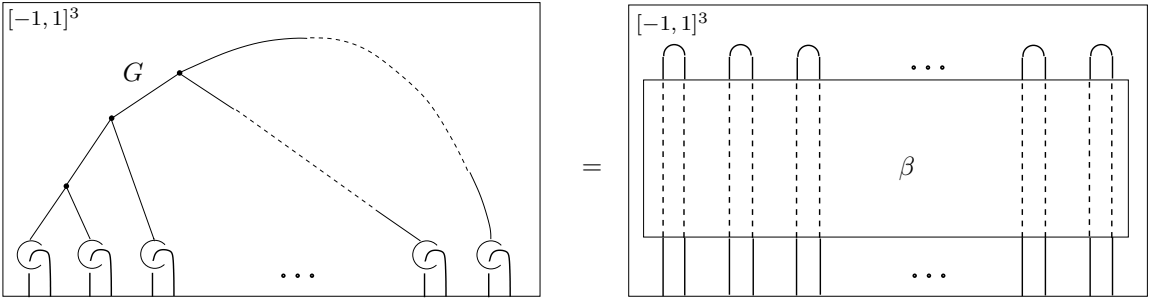}
\adjustrelabel <-0.2cm,0.1cm> {G}{$G$}
\adjustrelabel <.cm,-0.15cm> {B}{\scriptsize{$[-1,1]^3$}}
\adjustrelabel <.cm,-0.15cm> {B'}{\scriptsize{$[-1,1]^3$}}
\adjustrelabel <.2cm,-0.1cm> {beta}{$\beta$}
\adjustrelabel <.cm,-0.2cm> {=}{$=$}
\endrelabelbox}
\caption{The effect of the surgery along a tree clasper $G$.}
\label{fig:picture}
\end{figure}
More precisely, $\beta$ defines an element of the $(k+1)$-st term of the lower central series 
of the pure braid group $P_{2(k+2)}$ on $2(k+2)$ strands, 
and $\beta$ is given the Seifert framing of $[-1,1]^3$.
Now, let $\Sigma \subset M$ be the sub-surface of $\partial \hbox{N}(G)$ corresponding 
to $\left([-1,1]^2\times (-1)\right) \cap \left([-1,1]^3 \setminus\int\ \hbox{N}(\tau)\right)$. 
This is a disk with $2(k+2)$ holes and whose mapping class group $\mathcal{M}(\Sigma)$ can be identified
with $P_{2(k+2)} \times \Z^{2(k+2)}$, i.e$.$ the pure braid group on $2(k+2)$ framed strands. 
Thus, $P_{2(k+2)}$ is embedded as a subgroup of $\mathcal{M}(\Sigma)$ and, in particular, 
the braid $\beta$ defines an element of $\mathcal{M}(\Sigma)$. Then,
Figure \ref{fig:picture} can be interpreted as follows: $M_G$ is obtained from $M$ 
by cutting open along $\Sigma$ and re-gluing with $\beta$. 
Now, let $H$ be the regular neighborhood of $\Sigma$ in $M$: This is a genus $2(k+2)$ handlebody
which we identify with $\Sigma \times [-1,1]$.
Thus, $h:= (\beta \times (-1)) \cup (\Id \times 1)$ defines 
an orientation-preserving homeomorphism $\partial H \to \partial H$ 
such that $M_h = (M \setminus \int\ H) \cup_h H$ is homeomorphic to $M_G$.
Let also $\Sigma'$ be the surface $\partial H = \partial (\Sigma \times [-1,1])$ 
deprived of a small disk disjoint from $\Sigma \times (-1)$. 
Then, the inclusion of surfaces $\Sigma \subset \Sigma'$ 
induces a monomorphism $\mathcal{M}(\Sigma) \subset \mathcal{M}(\Sigma')$ at the level of mapping class groups, 
which has been studied first  by Oda \cite{Oda}. 
According to \cite{GH}, this sends $\gamma_{n+1}(P_{2(k+2)})$ to 
$\gamma_n(\mathcal{I}(\Sigma'))$ for all $n\geq 1$.
We conclude that $M \leadsto M_h$ is a Torelli surgery with $h \in \gamma_k(\mathcal{I}(\partial H))$.
\end{proof}

\bibliographystyle{amsalpha}

\end{document}